\documentclass[journal]{IEEEtran}  

\usepackage{ifpdf}

\usepackage{cite}

\ifCLASSINFOpdf
   \usepackage[pdftex]{graphicx}
   \DeclareGraphicsExtensions{.pdf,.jpeg,.png}
\else
   \usepackage[dvips]{graphicx}
   \DeclareGraphicsExtensions{.eps}
\fi
\usepackage{amsmath, amssymb, bm}
\usepackage{array}

\ifCLASSOPTIONcompsoc
  \usepackage[caption=false,font=normalsize,labelfont=sf,textfont=sf]{subfig}
\else
  \usepackage[caption=false,font=footnotesize]{subfig}
\fi

\usepackage{stfloats}
\ifCLASSOPTIONcaptionsoff
  \usepackage[nomarkers]{endfloat}
 \let\MYoriglatexcaption\caption
 \renewcommand{\caption}[2][\relax]{\MYoriglatexcaption[#2]{#2}}
\fi
\usepackage{url}

\usepackage{float}
\usepackage{textcomp}

\newtheorem{theorem}{Theorem}
\newtheorem{lemma}{Lemma}
\newtheorem{definition}{Definition}
\newtheorem{remark}{Remark}
\newtheorem{example}{Example}

\begin{document}

\makeatletter
\newcommand{\norm}[1]{\left\lVert#1\right\rVert}
\makeatother
%
\title{A note on uniform exponential stability of linear periodic time-varying systems}
%
%
%

\author{Robert~Vrabel
\thanks{\textcopyright 2019 IEEE.  Personal use of this material is permitted. Permission from IEEE must be obtained for all other uses, in any current or future media, including reprinting/republishing this material for advertising or promotional purposes, creating new collective works, for resale or redistribution to servers or lists, or reuse of any copyrighted component of this work in other works.}
\thanks{Digital Object Identifier 10.1109/TAC.2019.2927949}
}

%
%

\markboth{IEEE Transactions on Automatic Control }%
{Vrabel}
%



\maketitle

\begin{abstract}
In this paper we derive new criterion for uniform stability assessment of the linear periodic time-varying systems $\bm{\dot x=A(t)x,$\ $A(t+T)=A(t).}$ As a corollary, the lower and upper bounds for the Floquet characteristic exponents are established. The approach is based on the use of logarithmic norm of the system matrix $\bm{A(t).}$ Finally we analyze the robustness of the stability property under external disturbance. 
\end{abstract}

\begin{IEEEkeywords}
Linear periodic time-varying system, uniform exponential stability, uniform stability, logarithmic norm.
\end{IEEEkeywords}

%
\IEEEpeerreviewmaketitle


\section{Theory about general linear time-varying systems}
 
\subsection{Introduction}
Stability analysis for linear time-varying (LTV) systems  is of constant interest in the control community. One reason is the growing importance of adaptive controllers for which underlying closed-loop adaptive system is time-varying and linear \cite{Gao_et_al, Ioannou, Shamma}. The second one is that the LTV systems naturally arise when one linearizes nonlinear systems about a non constant nominal trajectory. In contrast the linear time-invariant (LTI) cases which have been thoroughly understood in the analysis and synthesis, many properties of the LTV systems are still not completely resolved. In this context, the system stability analysis can serve as an appropriate example. 

The stability characteristics of a linear time-invariant (LTI) system of ordinary differential equations $\dot x =Ax$ can be characterized completely by the placement of the eigenvalues of the system matrix $A.$
For systems described by
\begin{equation}\label{state_equation}
\dot x =A(t)x, \quad t\geq t_0\,(\geq 0),
\end{equation}
one would intuitively expect that  if,  for each  $t,$  the frozen-time system  is
stable of any kind, then the time-varying system should also be stable provided  $A(t)$ is bounded. However
these conditions are still not strong enough to guarantee the uniform exponential stability (Example~\ref{example_contra}) and additional restrictions suitably constraining the rate of variation in $A(t)$ have to be imposed.  The best known results were given by C.~A.~Desoer \cite{Desoer1}, W.~A.~Coppel \cite{Coppel1}  and H.~H.~Rosenbrock \cite{Rosenbrock} in their studies of slowly varying systems. The results are summarized and slightly strengthened in \cite[Theorem~3.2]{Ilchmann}.  For illustration purpose, in the following theorem we present two criteria; for some other frozen-time methods for LTV systems see also e.~g.\cite{Izobov}.
\begin{theorem}
Suppose that $A(t)$ is (piecewise) continuous matrix function $A(\cdot):[0,\infty)\to \mathbb{R}^{n\times n}$ 
which satisfies:
\begin{itemize}
\item[i)]there exists $M>0$ such that $\norm{A(t)}<M$ for all $t\geq0,$
\item[ii)] there exists $\alpha>0$ such that the spectrum 

$
\sigma_{A(t)}\subset\left\{z\in\mathbb{C};\, \Re\{z\}<-\alpha \right\}$ for all $t\geq 0.
$
\end{itemize}
Then any of the following conditions guarantees uniform exponential stability of (\ref{state_equation}):
\begin{itemize}
\item[C1)] $\alpha>4M$ for all $t\geq0;$
\item[C2)] $A(\cdot)$ is piecewise differentiable and 
\[
\norm{\dot A(t)}<\frac{2}{2n-1}\frac{\alpha^{4n-2}}{2M^{4n-4}}\ \mathrm{for\ all}\ t\geq0.
\]
\end{itemize}
\end{theorem}
\begin{example}\label{example_contra}
For the linear state equation (\ref{state_equation}) with
\begin{equation*}
A_\beta(t)=\left[
\begin{array}{cc} 
-1+\beta\cos^2t & 1-\beta\sin t\cos t \\ 
-1-\beta\sin t\cos t & -1+\beta\sin^2t
\end{array}
\right],\ t\geq0,
\end{equation*}
where $\beta$ is a positive constant, the pointwise eigenvalues are constants, given by 
$
\sigma_{A_\beta(t)}=\frac{1}{2}\left({\beta-2\pm\sqrt{\beta^2-4}}\right).
$ It is not difficult to verify that the fundamental matrix
\[
\Phi_\beta(t,0)= \left[
\begin{array}{cc} 
e^{(\beta-1)t}\cos t & e^{-t}\sin t \\ 
-e^{(\beta-1)t}\sin t  & e^{-t}\cos t
\end{array}
\right].
\]
Thus while the pointwise eigenvalues of $A_\beta(t)$ have negative real parts if $0 < \beta < 2,$ the
state equation has unbounded solutions if $\beta>1,$ \cite{Rugh}. Now if we set $\beta=1.5$ and $t=t^*=\pi/2,$ we get
$\norm{A_\beta(t^*)}_2=1.7808,$ $\norm{\dot A_\beta(t^*)}_2=0.5000,$
$\sigma_{A_\beta(t^*)}=\{ -0.2500 + 0.6614i, \, -0.2500 - 0.6614i\},$
and so the system does not satisfy neither sufficient condition C1 nor C2. In general, it is difficult to specify exact upper bounds on $\norm{A(t)}$ and $\norm{\dot A(t)}$ for all $t\geq0.$

Principally different approach to study of stability of LTV systems is based on the analysis of the small perturbation of the stable nominal system $\dot x=A_{\mathrm{nom}}x.$ As is shown in \cite{Bellman} and \cite{Coppel2}  the perturbed system $\dot x=[A_{\mathrm{nom}}+B(t)]x$ preserves the uniform and uniform exponential stability if $\int\limits_0^\infty\norm{B(\tau)}d\tau<\infty.$
\end{example}
\subsection{Notation and definitions}
\begin{definition}[\cite{Brockett1}, \cite{Coppel2}]\label{definition_stability}
Let $X(t)$ is a fundamental matrix solution for (\ref{state_equation}) and $\Phi(t,\tau)\triangleq X(t)X^{-1}(\tau)$ denotes its corresponding state-transition matrix. Then the system (\ref{state_equation}) is 
\begin{itemize}
\item[(US)] uniformly stable if and only if there exists a positive constant $K$ such that 
\begin{equation*}
\norm{\Phi(t,\tau)}\leq K\ \mathrm{for}\ t_0\leq\tau\leq t<\infty,
\end{equation*}
\item[(UES)] uniformly exponentially stable if and only there exist positive constants
$K$, $\tilde\alpha$ such that
\begin{equation*}
\norm{\Phi(t,\tau)}\leq Ke^{-\tilde\alpha(t-\tau)}\ \mathrm{for}\ t_0\leq\tau\leq t<\infty.
\end{equation*}
\end{itemize}
\end{definition}
We will derive results for unspecified vector norm on $\mathbb{R}^n$, $\norm{\cdot}.$ For the matrices, as an operator norm is used the induced norm, $\norm{A}=\max\limits_{\norm{x}=1}\norm{Ax}$.  We will use for both vector norm and matrix operator norm the same notation but it will always be clear from the context that norm is being used. In particular cases we will consider the three most common vector norm - $\norm{x}_1,$ $\norm{x}_2$ and $\norm{x}_{\infty}.$ We denote by $\mu[A(t)],$ $t\geq0,$  the logarithmic norm (LN) of a continuous matrix function $A(t),$ $A(\cdot):[0,\infty)\to \mathbb{R}^{n\times n}$ defined as
\[
\mu[A(t)]\triangleq\lim\limits_{h\to 0^+}\frac{\norm{I_n+hA(t)}-1}{h},
\]
where $I_n$ denotes the identity on $\mathbb{R}^n$ (see Table~\ref{table:norms}). We note here that the LN $\mu$ is not a norm in the usual sense, because it can take negative values.
\begin{table}[H]
\renewcommand{\arraystretch}{1.3}
\caption{Logarithmic norms for the vector norms $\norm{\cdot}_1,$ $\norm{\cdot}_2$ and $\norm{\cdot}_\infty,$ \cite[p.~54]{Afanasiev}, \cite[p.~33]{Desoer_Vidyasagar}.}
\label{table:norms}
\centering
\begin{tabular}{c c}
\hline
\bfseries Vector norm &  \bfseries Logarithmic norm  \\
\hline\hline
$\norm{x}_1=\sum\limits_{i=1}^{n}|x_i|$ & $\mu_1[A]=\max\limits_{1 \leq j \leq n}\left(a_{jj}+\sum\limits_{i\neq j}|a_{ij}|\right)$   \\
 &  \\
$\norm{x}_2=\sqrt{\sum\limits_{i=1}^{n}x_i^2}$ & $\mu_2[A]=\frac12\lambda_{\max}\left({A+A^T}\right)$   \\
 &  \\
$\norm{x}_{\infty}=\max\limits_{1 \leq i \leq n}{|x_i|}$ & $\mu_{\infty}[A]=\max\limits_{1 \leq i \leq n}\left(a_{ii}+\sum\limits_{j\neq i}|a_{ij}|\right)$   \\  
\hline
\end{tabular} 
\end{table}
In Table~\ref{table:norms} and elsewhere in the paper, the superscript 'T' denotes transposition,  the number $\lambda_{\max}(A+A^T)$ is the maximum eigenvalue of the matrix $A+A^T.$
\begin{remark}\label{LTI_lognorm}
Note that the value $\mu[A]$ may depends on the used vector norm, see an example in \cite[p.~56]{Afanasiev}. Thus, we can verify whether the LTI system $\dot x=Ax$ is stable or not by means of the vector norm with negative value of $\mu[A],$ see Lemma~\ref{lognorm_properties} (P3) below. For a Hurwitz matrix $A$ we obtain such LN for a vector norm $\norm{x}_H\triangleq\sqrt{x^THx},$ where the symmetric positive definite matrix $H$ satisfies the Lyapunov equation $A^TH + HA = -2I_n.$ Then $\mu_H[A]=-1/\lambda_{\max}(H),$  see Lemma~2.3 in \cite{Hu_Liu}. Thus, the stability in terms of LN becomes a topological notion, while the spectrum $\sigma_A=\{\lambda\in\mathbb{C}:\, \lambda\ \mathrm{is\ an\ eigenvalue\ of}\ A\}$ is topologically invariant.
\end{remark}
Now we summarize the important properties of the LN useful for the stability analysis of linear dynamical systems.
\begin{lemma}[\cite{Coppel2, Desoer2, Desoer_Vidyasagar, Soderlind1, Soderlind2}]\label{lognorm_properties}
\begin{itemize}
\item[P1)] $-\mu[-A]\leq\mu[A];$ $\vert\mu[A]-\mu[B]\vert\leq\norm{A-B}$ for any given $n\times n$ matrices $A$ and $B;$
\item[P2)] Let $X(t),$ $t\geq0$ is a fundamental matrix solution for $\dot x=A(t)x.$ Then 
\begin{equation}\label{property_P2}
e^{-\int\limits_{\tau}^t \mu[-A(s)]ds}\leq\norm{X(t)X^{-1}(\tau)}\leq e^{\int\limits_{\tau}^t \mu[A(s)]ds}
\end{equation}
for all $t_0\leq\tau\leq t<\infty;$
\item[P3)] The solution of (\ref{state_equation}) satisfies for all $t\geq t_0$ the inequalities
\[
\norm{x(t_0)}e^{-\int\limits_{t_0}^t \mu[-A(s)]ds}\leq\norm{x(t)}\leq\norm{x(t_0)}e^{\int\limits_{t_0}^t \mu[A(s)]ds}.
\]
\end{itemize}
\end{lemma} 
\section{Theory about linear periodic time-varying systems}
Although precise stability assessment for general LTV systems is very difficult, the stability of linear periodic time-varying (LPTV) systems 
\begin{equation}\label{state_equation_per}
\dot x=A(t)x,\quad A(t+T)=A(t)\ \mathrm{for\ some}\ T>0, 
\end{equation}
can be determined using the Floquet theory, which states that for every LPTV system the associated state-transition matrix can be expressed as 
\begin{equation}\label{state_transition_Lyap}
\Phi(t,\tau)=P^{-1}(t)e^{R(t-\tau)}P(\tau),\  0\leq\tau\leq t<\infty 
\end{equation}
via a stability preserving Lyapunov transformation $z=P(t)x,$ where $P(t)$ is continuously differentiable, nonsingular and periodic matrix function which reduces the stability analysis to an analysis of the LTI system $\dot z=Rz,$ where we define the $n\times n$ constant matrix $R$ by setting $e^{RT}=\Phi(T,0).$ Consequently, the stability of the original LPTV system is equivalent to that of the LTI system  $\dot z=Rz.$ The eigenvalues of $R$ are known as the Floquet characteristic exponents (FCEs) \cite{Brockett1}, \cite{Chicone}, \cite{Coppel2}. Unfortunately, the application of this theory is hampered by the fact that, in general, the FCEs is hard to determine, namely, the Lyapunov transformation $P(t)$ is defined via $P^{-1}(t)=\Phi(t,0)e^{-Rt}.$  

A comprehensive Floquet theory including Lyapunov transformations was developed and their various stability preserving properties were analyzed in \cite{DaCunha}. Colaneri \cite{Colaneri} addresses a few theoretical aspects of LPTV systems and methodology which can be useful to characterize and extend other concepts usually exploited in the time-invariant case only. On the other hand, relating computational and numerical aspects, in \cite{Tian_Wang}, the FCEs are directly calculated for the special types of system matrices, when the coefficient matrices are triangular. In \cite{Garcia}, based on the solution of linear differential Lyapunov matrix equation, necessary and sufficient numerical conditions for asymptotic stability of LTV systems are given. 

Our aim in the present paper is to provide a conceptually new approach to study of general LPTV systems allowing to estimate the norm of state-transition matrix and subsequently its stability property without knowing the fundamental matrix solution, purely on the basis of the matrix $A(t)$ entries (Theorem~\ref{main_theorem}). This approach has the advantage of avoiding the need of calculation of Lyapunov transformation. Moreover, the developed technique allows us to find the upper and lower bounds for the solutions of LPTV systems (\ref{state_equation_per}) in Lemma~\ref{main_lemma} and for the FCEs in Remark~\ref{FCE_bounds}. As is shown, the accuracy of the achieved estimates depends on the used vector norm in $\mathbb{R}^n.$ Because the spectrum of a matrix is invariant to the change of norm on $\mathbb{R}^n,$ this problem can be formulated as an optimization problem of finding vector norm on $\mathbb{R}^n$ minimizing (separately) $\lambda^-$ and $\lambda^+.$ Definitions of these and other important constants and concepts are given in the following subsection.
\subsection{Notation (continued)}
Let
\[
\Pi^+(t)\triangleq\int\limits_{t_0}^{t}\mu[A(s)]ds, \quad \Pi^-(t)\triangleq\int\limits_{t_0}^{t}\mu[-A(s)]ds,
\]
\[
\lambda^+\triangleq\Pi^+(t_0+T)/T,\quad \lambda^-\triangleq\Pi^-(t_0+T)/T
\]
\[
\delta^*_U\triangleq\min\left\{\delta:\,\Pi^*(t)\leq\lambda^*(t-t_0)+\delta, \forall t\in[t_0,t_0+T] \right\}, 
\] 
\[
\delta^*_L\triangleq\max\left\{\delta:\, \Pi^*(t)\geq\lambda^*(t-t_0)+\delta, \forall t\in[t_0,t_0+T] \right\}, 
\]
$ *=+,-.$ All functions and constants are well-defined because $\mu[A(t)]$ is continuous as follows from Lemma~\ref{lognorm_properties} (P1) and from the assumption of continuity of matrix function $t\to A(t).$  The functions $e^{\lambda^*(t-t_0)+\delta^*_L}$ and  $e^{\lambda^*(t-t_0)+\delta^*_U},$ $ *=+,-$ will be called a lower and upper barrier function, respectively. The constants $\delta^*_L,$ $\delta^*_U$ can be calculated by applying global extrema--searching procedure for the function $\Pi^*(t)-\lambda^*(t-t_0)$ on the interval $[t_0,t_0+T].$
\subsection{Auxiliary results}
In all lemmas below it is assumed that $A(t)$ is periodic with period $T>0.$
\begin{lemma}
$
\ \delta^*_U\geq0,\quad \delta^*_L\leq0,\quad *=+,\,-.
$
\end{lemma}
\begin{IEEEproof}
Both inequalities follows immediately from the fact that $\Pi^*(t_0+T)=\lambda^*T,$ $*=+,\,-.$
\end{IEEEproof}
\begin{lemma}\label{pi_inequalities1}
$-\Pi^-(t)\leq\Pi^+(t)$ for all $t\geq t_0.$
\end{lemma}
\begin{IEEEproof}
This property follows from Lemma~\ref{lognorm_properties} (P1).
\end{IEEEproof}
\begin{lemma}\label{pi_inequalities2}
For all $t\geq t_0$ is $\Pi^*(t)\leq\lambda^*(t-t_0)+\delta^*_U$ and $\Pi^*(t)\geq\lambda^*(t-t_0)+\delta^*_L,$ $*=+,\,-.$
\end{lemma}
\begin{IEEEproof}
Let $\hat t\in[t_0,\infty)$ is chosen arbitrarily. Then $\hat t\in[t_0+(k-1)T, t_0+kT)$ for some $k\geq1.$ As follows from the definition of LN, $\mu[A(t)]$ is also periodic which yields 
\[
\Pi^*(\hat t)= (k-1)\Pi^*(t_0+T)+\int\limits_{t_0}^{\hat t-(k-1)T}\mu[*A(s)]ds
\]
\[
=(k-1)\Pi^*(t_0+T)+\Pi^*(\hat t-(k-1)T),\ *=+,\,-.
\]
Now, because $\hat t-(k-1)T\in[t_0,t_0+T),$ we have that
\[
\Pi^*(\hat t)\leq (k-1)\Pi^*(t_0+T)+\lambda^*(\hat t-(k-1)T-t_0)+\delta^*_U
\]
\[
=(k-1)\Pi^*(t_0+T)+\frac{\Pi^*(t_0+T)}{T}(\hat t-(k-1)T-t_0)+\delta^*_U
\]
$=\lambda^*(\hat t-t_0)+\delta^*_U,$ what we had to prove. The extension of the inequality from the interval $[t_0,t_0+T]$ on the whole time interval $[t_0,\infty)$  for the lower bounds of $\Pi^*$ can be proved in a similar manner.
\end{IEEEproof}
\begin{lemma}\label{main_lemma}
Let us consider (\ref{state_equation_per}) with an initial state $x(t_0)\in\mathbb{R}^n.$ Then
\[
\norm{x(t_0)}e^{-\lambda^-(t-t_0)-\delta^-_U}\leq \norm{x(t_0)}e^{-\Pi^-(t)}\leq\norm{x(t)}
\]
\begin{equation}\label{estimates}
\leq\norm{x(t_0)}e^{\Pi^+(t)}\leq\norm{x(t_0)}e^{\lambda^+(t-t_0)+\delta^+_U}
\end{equation}
for all $t\geq t_0.$
\end{lemma}
\begin{IEEEproof} 
The claim of lemma follows from Lemma~\ref{lognorm_properties} (P3) and Lemma~\ref{pi_inequalities2}. Lemma~\ref{pi_inequalities1} guarantees that the inequalities in (\ref{estimates}) make sense. 
\end{IEEEproof}
\begin{remark}\label{FCE_bounds}
As a corollary we obtain for FCEs the inclusion 
\[
\sigma_{R}\subset\{z\in\mathbb{C}: -\lambda^-\leq\Re\{z\}\leq \lambda^+ \}.
\]
In fact, if there was an eigenvalue $\hat\lambda\in\sigma_R$ such that $\Re\{\hat\lambda\}>\lambda^+$ (analogously for $\Re\{\hat\lambda\}<-\lambda^-$), then, taking into account (\ref{state_transition_Lyap}), there would be a solution $\hat x(t)$ of (\ref{state_equation_per}) with $\norm{\hat x(t)}=\norm{\hat x(t_0)}o(e^{(\Re\{\hat\lambda\}+\varepsilon)(t-t_0)})$ as $t\to\infty$ for arbitrarily small constant $\varepsilon>0$ which contradicts with (\ref{estimates}). Here the asymptotics of $\hat x(t)$ is expressed by the "little-o"  Bachmann–-Landau notation.
\end{remark}
\subsection{Main results}
The sufficient conditions for stability of the LPTV systems can be expressed in terms of an integral over one period $T$ of the LN $\mu[A(t)]$ or $\mu[-A(t)].$
\begin{theorem}\label{main_theorem}
If for some vector norm and associated induced norm for matrices is 
\begin{itemize}
\item[I)] $\Pi^+(t_0+T)<0,$ then the LPTV system (\ref{state_equation_per}) is UES;
\item[II)] $\Pi^+(t_0+T)=0,$ then the LPTV system (\ref{state_equation_per}) is US;
\item[III)] $\Pi^-(t_0+T)<0,$ then the LPTV system (\ref{state_equation_per}) is unstable (specifically, the norms of all nonzero solutions converge to infinity as $t\to\infty$).
\end{itemize}
\end{theorem}
\begin{IEEEproof} 
\underline{I)+II)}:\ For all $t_0\leq\tau\leq t<\infty$ we have, by (\ref{property_P2}),
\begin{equation*}
\norm{X(t)X^{-1}(\tau)}\leq e^{\int\limits_{\tau}^t \mu[A(s)]ds}=e^{\int\limits_{t_0}^t\mu[A(s)]ds-\int\limits_{t_0}^\tau\mu[A(s)]ds}
\end{equation*}
\[
=e^{\Pi^+(t)-\Pi^+(\tau)}\leq e^{(\lambda^+(t-t_0)+\delta^+_U)-(\lambda^+(\tau-t_0)+\delta^+_L)}
\]
\begin{equation*}
=e^{\lambda^+(t-\tau)}e^{\delta^+_U-\delta^+_L}=e^{\frac{\Pi^+(t_0+T)}{T}(t-\tau)}e^{\delta^+_U-\delta^+_L},
\end{equation*}
that is, we set $K=e^{\delta^+_U-\delta^+_L}$ and $\tilde\alpha=-\frac{\Pi^+(t_0+T)}{T}$ in Definition~\ref{definition_stability}.

\underline{III)}:\ From the left inequalities in (\ref{estimates}) and because $\lambda^-=\Pi^-(t_0+T)/T$ by definition, it follows that the norm of each nonzero solution of (\ref{state_equation_per}) converges to infinity as $t\to\infty,$ or alternatively,
analogously as above,
\begin{equation*}
\norm{X(t)X^{-1}(\tau)}\geq e^{\Pi^-(\tau)-\Pi^-(t)}
\end{equation*}
\[
\geq e^{-\frac{\Pi^-(t_0+T)}{T}(t-\tau)}e^{\delta^-_L-\delta^-_U}\to\infty
\]
as $t\to\infty$ for every fixed $\tau\geq t_0.$
\end{IEEEproof}
\begin{remark} 
Combining Lemma~\ref{lognorm_properties} (P2) with \cite[Lemma~2 (Item~3)]{Zhou} and \cite[Lemma~5]{Zhou} for $\mu(t)\triangleq\mu[A(t)]$ we get another justification of the sufficient condition for uniform exponential stability in Part~I of theorem with the difference that we have also derived the values $K$ and $\tilde\alpha$ from Definition~\ref{definition_stability} ($e^{\beta}$ and $\alpha$ in \cite{Zhou}) and which are generally classified as "difficult to obtain". 
\end{remark} 
\begin{remark} The connection of Theorem~\ref{main_theorem} with LTI systems $\dot x=Ax,$ which can be considered as LPTV systems with any period $T>0:$
\begin{itemize}
\item[I)] Remark~\ref{LTI_lognorm} implies for a Hurwitz matrix $A$ 
\[
\Pi^+(t_0+T)=\int\limits_{t_0}^{t_0+T}\mu_{H}[A]ds=-T/\lambda_{\max}(H)<0
\]
This shows that the condition $\Pi^+(t_0+T)<0$ is at the same time also necessary condition to be the LTI system $\dot x=Ax$ UES.
\item[II)] Let a nonsingular $n\times n$ real matrix $P$ and matrix $A$ are such that $J_A\triangleq PAP^{-1}$ has a Jordan normal form, where each block is of the form 
$
J_1=\left[
\begin{array}{cc} 
a & b \\
 -b & a   
\end{array}
\right]
$
or $J_2=[\tilde\lambda]$ with $a\leq0,$ $b\neq0$ and $\tilde\lambda\leq0.$ Define on $\mathbb{R}^n$ the $\norm{x}_{P^TP} (=\norm{Px}_2).$ Then, from the equality $\mu_{P^TP}[A]=\mu_2[J_A]$ \cite{Desoer2} we get $\mu_{P^TP}[A]=0$ for block diagonal matrix $J_A$ with at least one $J_1$ or $J_2$ with $a=0$ or $\tilde\lambda=0,$ that is,
\[
\Pi^+(t_0+T)=\int\limits_{t_0}^{t_0+T}\mu_{P^TP}[A]ds=0;
\]
Thus, the condition $\Pi^+(t_0+T)=0$ is also necessary for the systems with Jordan normal form described above and which are US but not UES.
\item[III)] If all eigenvalues of $A$ have positive real part then $(-A)$ is a Hurwitz matrix because $\sigma_{-A}=-\sigma_A$ as follows from the equality $\det(A-\lambda I_n)=(-1)^n\det((-A)-(-\lambda)I_n)$ \cite[p.~524]{Harville}. Remark~\ref{LTI_lognorm} yields that 
\[
\Pi^-(t_0+T)=\int\limits_{t_0}^{t_0+T}\mu_{\tilde H}[-A]ds=-T/\lambda_{\max}(\tilde H)<0.
\]
Thus, $\Pi^-(t_0+T)<0$ establishes also a necessary condition for instability of LTI systems $\dot x=Ax,$ with $(-A)$ being a Hurwitz matrix.
\end{itemize}
\end{remark}
Revisiting Example~\ref{example_contra} in the light of Theorem~\ref{main_theorem} we see 
\[
\Pi^+_{\beta}(2\pi)=\int\limits_0^{2\pi}\mu_2[A_{\beta}(s)]ds=\int\limits_0^{2\pi}\max\sigma_{\frac12(A_\beta(s)+A^T_\beta(s))}ds
\]
\[
=2\pi\max\{-1,\beta-1\},
\]
\[
\Pi^-_{\beta}(2\pi)=\int\limits_0^{2\pi}\mu_2[-A_{\beta}(s)]ds=\int\limits_0^{2\pi}\max\sigma_{\frac12(-A_\beta(s)-A^T_\beta(s))}ds
\]
\[
=2\pi\max\{1,1-\beta\},
\]
and so $\Pi^+_{\beta}(2\pi)<0$ if $\beta<1$ (UES system), $\Pi^+_{\beta}(2\pi)=0$ if $\beta=1$ (US system). The sufficient condition for instability is not fulfilled because there is also exponentially stable mode in the system, not influenced by the parameter $\beta.$
\begin{example}
As an illustrative example let us consider for $t\geq0$ the LPTV system with
\begin{equation}\label{eq:example2} 
\arraycolsep=1.0pt\def\arraystretch{1.2}
A(t)=\left[
\begin{array}{cc} 
-11/2+(15/2)\sin12t & (15/2)\cos12t \\
 (15/2)\cos12t & -41/2-(15/2)\sin12t   
\end{array}
\right]
\end{equation}
For comparison purpose we calculate the barrier functions for the two vector norms, $\norm{\cdot}_1$ and $\norm{\cdot}_2.$ Meaning the lower and upper barrier function is obvious from Fig.~\ref{functions_mu1_approx}. We have that
\[
\mu_1[A(t)]=-11/2+(15/2)(\sin12t+\left\vert\cos12t\right\vert)\ (=\mu_\infty[A(t)])
\]
\[
\lambda^+=\frac{6}{\pi}\int\limits_0^{\pi/6}\left[-11/2+(15/2)(\sin12s+\left\vert\cos12s\right\vert)\right]ds
\]
$=-0.7253$ and $\delta^+_U=1.2872.$ Analogously,
\[
\lambda^-=\frac{6}{\pi}\int\limits_0^{\pi/6}\left[41/2+(15/2)(\sin12s+\left\vert\cos12s\right\vert)\right]ds
\]
$=25.2747,$ and $\delta^-_U=1.2871.$ For the LN $\mu_2$ we analogously obtain:
\[
\mu_2[A(t)]=\frac{15\sqrt{2}\sqrt{\sin12t+1}}{2}-13\in[-13, 2],
\]
\[
\lambda^+=\frac{6}{\pi}\int\limits_0^{\pi/6}\left[\frac{15\,\sqrt{2}\,\sqrt{\sin 12s+1}}{2}-13\right]ds
\]
$=-3.4507,$ and $\delta^+_U= 0.9337;$
\[
\lambda^-=\frac{6}{\pi}\int\limits_0^{\pi/6}\left[\frac{15\,\sqrt{2}\,\sqrt{\sin 12s+1}}{2}+13\right]ds
\]
$=22.5493,$ and $\delta^-_U= 1.0441.$
\begin{figure}[!t]
\centering
\includegraphics[width=2.5in]{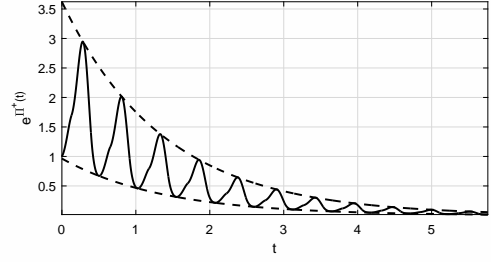}
\caption{The functions $e^{\Pi^+(t)}=e^{\int\limits_0^t \mu_1[A(s)]ds}$ with
$\mu_1[A(t)]=-11/2+(15/2)\big(\sin12t+\left\vert\cos12t\right\vert\big)$ (the solid line) and the lower and upper barrier functions $e^{\lambda^+t+\delta^+_L}=e^{-0.7253t-0.0364},$ $e^{\lambda^+t+\delta^+_U}=e^{-0.7253t+1.2872}$  (the dashed lines).}
\label{functions_mu1_approx}
\end{figure}
\begin{figure*}[ht]
\centering
\subfloat[Simulation for $\norm{\cdot}_1$]{\includegraphics[width=2.5in]{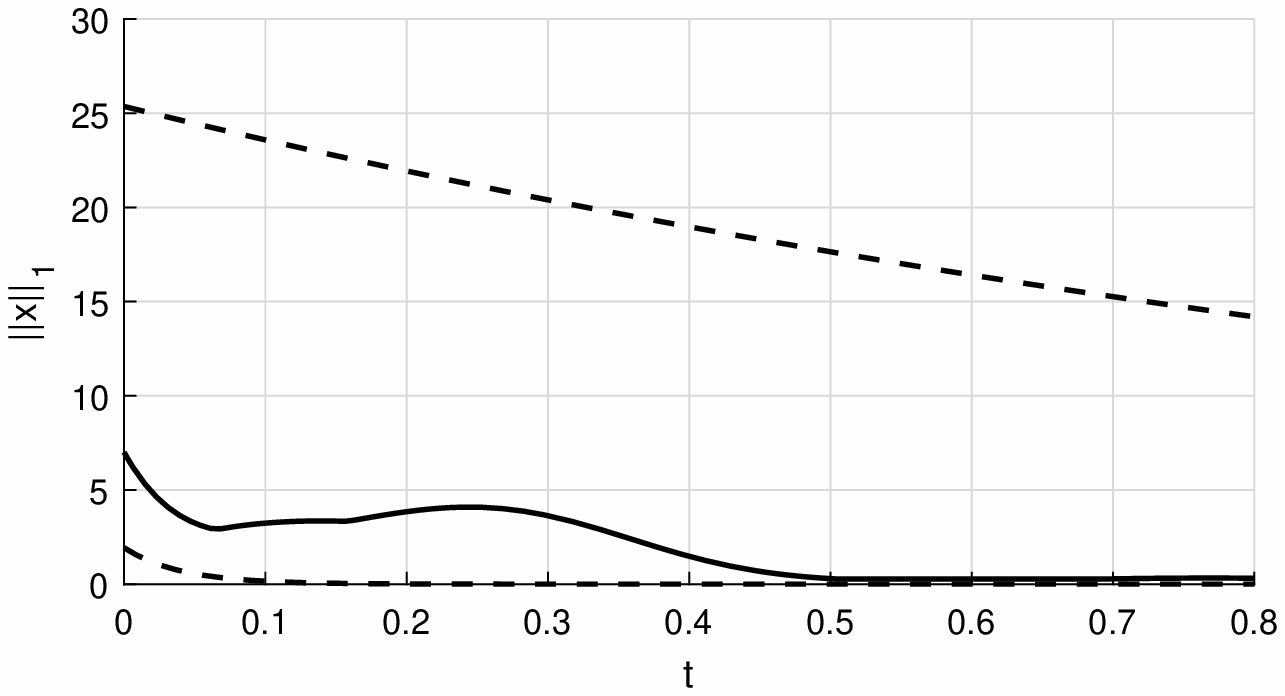}
}
\hfil
\subfloat[Simulation for $\norm{\cdot}_2$]{\includegraphics[width=2.5in]{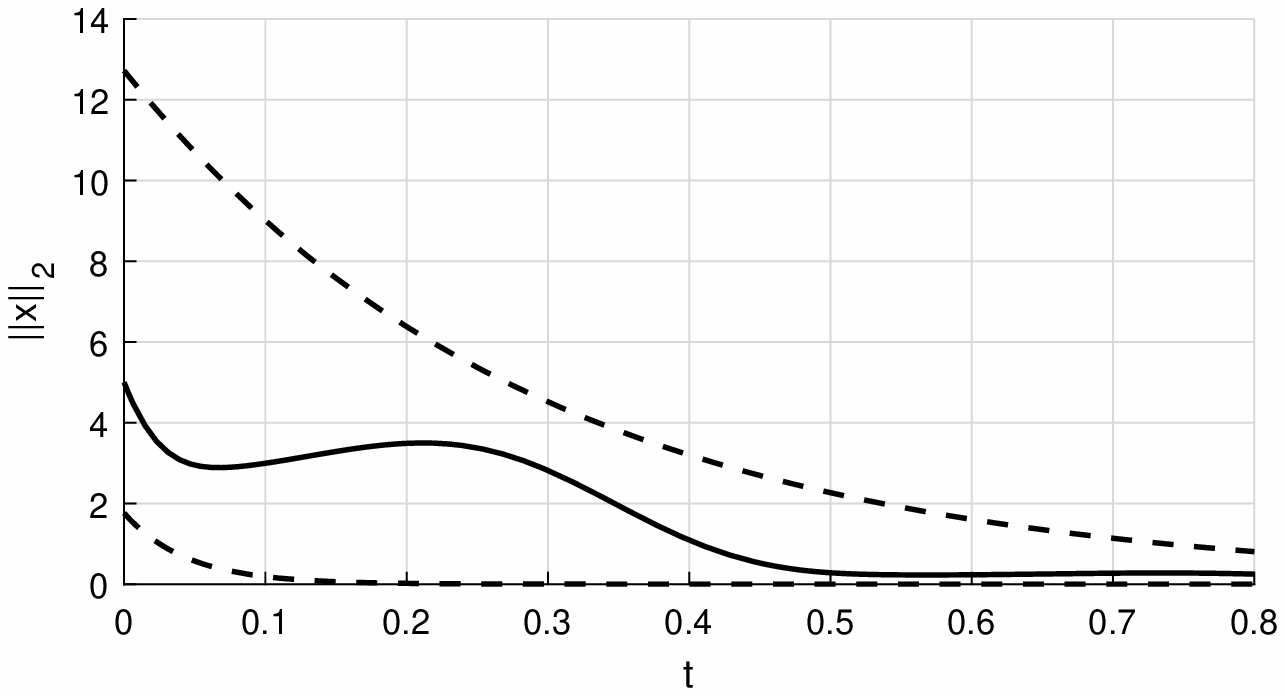}
}
\caption{The solution of (\ref{eq:example2}) with the initial state $x(0)=(-4,\ 3)^T$ (the solid line); the bounds for $\norm{x(t)}_1: (|-4|+|3|)e^{-0.7253t+1.2872}$, $(|-4|+|3|)e^{-25.2747t-1.2871}$ (the dashed lines); the bounds for $\norm{x(t)}_2: \sqrt{25}e^{-3.4507t+0.9337}$, $\sqrt{25}e^{-22.5493t-1.0441}$ (the dashed lines).}
\label{figsim}
\end{figure*}
\end{example}
In this example, the use of the Euclidean vector $\mu_2[A(t)]$ provides the better information regarding the position of the FCEs ($=\sigma_R$) in the complex-plane as those given by octahedral norm $\mu_1[A(t)],$ the vertical strip $\{z\in\mathbb{C}:\, -22.5493\leq\Re\{z\}\leq -3.4507 \}_{\mu_2}$ and $\{z\in\mathbb{C}:\, -25.2747\leq\Re\{z\}\leq -0.7253 \}_{\mu_1},$ respectively. The result of simulation is in Fig.~\ref{figsim}.
\subsection{Note regarding the robustness of exponentially stable LPTV systems against disturbances}
In this section we will analyze the stability properties of the LPTV systems affected by an external disturbance $d(t).$ Let us consider that the unperturbed system $\dot x=A(t)x$ is UES. What can we say about the asymptotic behavior of its  perturbation $\dot x=A(t)x+d(t)?$ This question represents one of the fundamental problems on the field of robust stability. The robustness of the systems' stability is not usually analyzed together with establishing the sufficient conditions ensuring the stability of some kind for the LTV systems. Among these include, e.~g., \cite{Ilchmann, Zhou, Safavi, Wang, Wu}.  We have the following result regarding asymptotic behavior of perturbed LPTV systems as $t\to\infty.$
\begin{theorem}
Let us consider the perturbed LPTV system,
\begin{equation}\label{state_equation_perturbed}
\dot x=A(t)x+d(t),\quad A(t+T)=A(t),\quad t\geq t_0\,(\geq0),
\end{equation}
where  $A(\cdot): [0,\infty)\to\mathbb{R}^{n\times n}$ and $d(\cdot): [0,\infty)\to\mathbb{R}^n$ are continuous and let 
\begin{itemize}
\item[A1)] $\Pi^+(t_0+T)<0,$ that is, the unperturbed system is UES, and 
\item[A2)]  $\norm{d(t)}\to 0$ as $t\to\infty.$
\end{itemize}
Then all solutions of (\ref{state_equation_perturbed}) converge to $0$ as $t\to\infty$ (not necessarily exponentially).
\end{theorem}
\begin{IEEEproof}
The solution of (\ref{state_equation_perturbed}) is given by the Lagrange's variation of constants formula,
\begin{equation*}
x(t)=X(t)X^{-1}(t_0)x(t_0)+X(t)\int\limits_{t_0}^t X^{-1}(\tau)d(\tau)d\tau
\end{equation*}
and so $\norm{x(t)}$
\[
\leq\norm{X(t)X^{-1}(t_0)}\norm{x(t_0)}+\int\limits_{t_0}^t \norm{X(t)X^{-1}(\tau)}\norm{d(\tau)}d\tau.
\]
Using the inequality $\norm{X(t)X^{-1}(\tau)}\leq e^{\frac{\Pi^+(t_0+T)}{T}(t-\tau)}e^{\delta^+_U-\delta^+_L}
$ from the proof of Theorem~\ref{main_theorem}, we have
\[
\norm{x(t)}\leq\norm{x(t_0)}e^{\frac{\Pi^+(t_0+T)}{T}(t-t_0)}e^{\delta^+_U}
\]
\[
+e^{\delta^+_U-\delta^+_L}\int\limits_{t_0}^t e^{\frac{\Pi^+(t_0+T)}{T}(t-\tau)}\norm{d(\tau)}d\tau
\]
\[
=\norm{x(t_0)}e^{\frac{\Pi^+(t_0+T)}{T}(t-t_0)}e^{\delta^+_U}
\]
\[
+e^{\delta^+_U-\delta^+_L}\frac{\int\limits_{t_0}^t e^{-\frac{\Pi^+(t_0+T)}{T}\tau}\norm{d(\tau)}d\tau}{e^{-\frac{\Pi^+(t_0+T)}{T}t}}.
\]
Applying the L'Hospital rule to the second term we get
\[
\lim\limits_{t\to\infty}\frac{\int\limits_{t_0}^t e^{-\frac{\Pi^+(t_0+T)}{T}\tau}\norm{d(\tau)}d\tau}{e^{-\frac{\Pi^+(t_0+T)}{T}t}}=\lim\limits_{t\to\infty}\frac{-\norm{d(t)}}{\Pi^+(t_0+T)/T}=0,
\] 
and so $\norm{x(t)}\to 0$ as $t\to\infty.$
\end{IEEEproof}
\begin{remark}
The perturbation $d(t)$ could have been replaced by a perturbation $\tilde d(x,t)$ which satisfies $\norm{\tilde d(x,t)}\leq \norm{d(t)}$ for all $x\in\mathbb{R}^n.$ Since this would have introduced no new ideas, we chose to present the notationally simpler case. 
\end{remark} 
The class of allowable disturbances of the form $d(t)$ preserving the convergence to $0$ of the solutions for the UES LPTV unperturbed systems is a little wider \cite{Strauss_Yorke} and contains also the functions that do not vanish at infinity.
\begin{theorem} 
Let the unperturbed system is UES. Then all solutions of perturbed LPTV system $\dot x=A(t)x+d(t),$ $A(t+T)=A(t)$ and perturbed LTI system $\dot x=Ax+d(t)$ converge to $0$ as $t\to\infty$ if and only if
\begin{equation}\label{diminishing}
\sup\limits_{0\leq\eta \leq1}\norm{\int\limits_t^{t+\eta}d(\tau)d\tau}\to 0\quad \mathrm{as} \quad t\to\infty. 
\end{equation}
\end{theorem} 
By \cite[Corollary~4.6]{Strauss_Yorke}, the class of allowable perturbations $d(x,t)$ contains also the functions of the form $D(t)k(x),$ where $D(t)$ is an $n\times n$ bounded matrix on $[0,\infty)$ whose columns satisfy (\ref{diminishing}) and $k: \mathbb{R}^n\to\mathbb{R}^n$  is continuous.
\section{Conclusion}
We derived, in Theorem~\ref{main_theorem}, the new criterion for uniform and uniform exponential stability of the linear periodic time-varying (LPTV) systems $\dot x=A(t)x,$ $t\geq t_0,$ $A(t+T)=A(t),$ $T>0$ without direct computing of the Floquet characteristic exponents (FCEs).
We have shown that the FCEs lie in the vertical strip $\{z\in\mathbb{C}: -\lambda^-\leq\Re\{z\}\leq \lambda^+ \}$ of the complex-plane. Here $\lambda^+=\frac1T\int\limits_{t_0}^{t_0+T}\mu[A(s)]ds$ and $\lambda^-=\frac1T\int\limits_{t_0}^{t_0+T}\mu[-A(s)]ds,$ where $\mu[\cdot]$ denotes the logarithmic norm (LN) of  matrix associated with an appropriately chosen vector norm. We also briefly discussed the persistence of the stability properties for the perturbed LPTV systems under the external disturbances $d(t).$ 
The fundamental advantage of the approach based on the use of LN is the fact that to estimate the norm of state-transition matrix $\Phi(t,\tau)$ for system $\dot x=A(t)x$  we do not need to know the fundamental matrix solution and all necessary estimates are based purely on the matrix $A(t)$ entries.


%


\ifCLASSOPTIONcaptionsoff
  \newpage
\fi

\end{document}